\documentclass[12pt]{amsart}%

\usepackage{amsmath}
\usepackage{amsfonts}
\usepackage{amssymb}
\usepackage{amstext}
\usepackage{graphicx}%
\providecommand{\U}[1]{\protect\rule{.1in}{.1in}}
\newtheorem{conjecture}{Conjecture}[section]

\newtheorem{theorem}[conjecture]{Theorem}

\newtheorem{corollary}[conjecture]{Corollary}

\newtheorem{proposition}[conjecture]{Proposition}

\numberwithin{equation}{section}

\begin{document}
\title{A generalization of almost Schur lemma on CR manifolds}
\author{Jui-Tang Chen$^{1\ast}$}
\address{$^{1}$Department of Mathematics, National Taiwan Normal University, Taipei
11677, Taiwan.}
\email{jtchen@ntnu.edu.tw}
\author{Nguyen Thac Dung$^{2}$}
\address{$^{2}$Department of Mathematics, National Taiwan University, Taipei 10617, Taiwam.}
\email{dungmath@yahoo.co.uk }
\author{Chin-Tung Wu$^{3\ast}$}
\address{$^{3}$Department of Applied Mathematics, National Pingtung University of
Education, Pingtung 90003, Taiwan.}
\email{ctwu@mail.npue.edu.tw }
\thanks{$^{\ast}$Research supported in part by NSC}

\begin{abstract}
In this paper, we study a general almost Schur Lemma on pseudo-Hermitian
(2n+1)-manifolds $(M,J,\theta)$ for $n\geq2.$ When the equality of almost
Schur inequality holds, we derive the contact form $\theta$ is pseudo-Einstein
and the pseudo-Hermitian scalar curvature is constant.

\end{abstract}

\maketitle

\section{Introduction}

In Riemannian manifolds, the classical Schur Lemma states that the scalar
curvature of an Einstein manifold of dimension $n\geq3$ must be constant. So
it is interesting to see the relation between scalar curvature and Einstein
condition. Recently, De Lellis and C. Topping \cite{lt} proved an almost Schur
Lemma assuming the nonnegative of Ricci curvature. Their result can be seen as
a quantitative version or a stability property of the Schur Lemma. Later, in
\cite{b}\cite{c1}\cite{c2} and \cite{gw}, the authors considered general
closed Riemannian manifolds, and obtained a generalization of the De
Lellis-Topping's theorem.

However, in the pseudo-Hermitian manifold, the pseudo-Einstein condition does
not imply the constant pseudo-Hermitian scalar curvature. This is because of
the appearance of torsion terms in the contracted Bianchi identity (\ref{6}).
Hence there is a natural question to ask under which condition a
pseudo-Einstein manifold has constant pseudo-Hermitian scalar curvature. More
general, how does the pseudo-Hermitian scalar curvature change when the
manifold is close to the pseudo-Einstein manifold. In \cite{csw}, the authors
addressed to this question and shown that if $Im \left(  A_{\alpha\beta
},^{\alpha\beta}\right)  =0$, where $A_{\alpha\beta}$ is the pseudo-Hermitian
torsion, then the answer is affirmative. In fact, the answer came from the
following CR almost Schur theorem in \cite{csw} on a closed pseudo-Hermitian
$(2n+1)$-manifold $M$ for $n\geq2$.

\begin{theorem}
\label{thm}(\cite{csw}) For $n\geq2,$ if $(M,J,\theta)$ is a closed
pseudo-Hermitian $(2n+1)$-manifold with $\operatorname{Im}\left(
A_{\alpha\beta},^{\alpha\beta}\right)  =0$ and%
\begin{equation}%
\begin{array}
[c]{c}%
(Ric-\frac{n+1}{2}Tor)(Z,Z)\geq0\text{ \ }\mathrm{for}\text{ }\mathrm{all}%
\text{ }Z\in T_{1,0}(M),
\end{array}
\label{01}%
\end{equation}
then%
\begin{equation}%
\begin{array}
[c]{c}%
\int_{M}(R-\overline{R})^{2}\leq\frac{2n(n+1)}{(n-1)\left(  n+2\right)  }%
\int_{M}%
%TCIMACRO{\tsum _{\alpha,\beta}}%
%BeginExpansion
{\textstyle\sum_{\alpha,\beta}}
%EndExpansion
|Ric_{\alpha\overline{\beta}}-\frac{R}{n}h_{\alpha\overline{\beta}}|^{2},
\end{array}
\label{1}%
\end{equation}
where $\overline{R}$ is the average value of the pseudo-Hermitian scalar
curvature $R$ over $M$. Moreover, equality holds then the contact form
$e^{\frac{1}{n+1}\varphi}\theta$ will be pseudo-Einstein.
\end{theorem}

In this paper, motivated by \cite{c1}, we are interested in a more general
curvature condition with respect to (\ref{01}). We prove a similar inequality
to (\ref{1}) with the inequality constant depending on the lower bound of
Webster Ricci tensor minus $\frac{n+1}{2}$ times torsion tensor and also on
the value of the first positive eigenvalue of the sub-Laplacian.

\begin{theorem}
\label{T}For $n\geq2,$ if $\left(  M,J,\theta\right)  $ is a closed
pseudo-Hermitian $(2n+1)$-manifold with%
\begin{equation}%
\begin{array}
[c]{c}%
(Ric-\frac{n+1}{2}Tor)(Z,Z)\geq-2K\left\vert Z\right\vert ^{2},
\end{array}
\label{0}%
\end{equation}
for all $Z\in T_{1,0}(M)$ and for some nonnegative constant $K$,\textrm{\ }%
then%
\begin{equation}%
\begin{array}
[c]{l}%
||R-\overline{R}||_{L^{2}}\leq\sqrt{\kappa}\left(  \int_{M}%
%TCIMACRO{\tsum _{\alpha,\beta}}%
%BeginExpansion
{\textstyle\sum_{\alpha,\beta}}
%EndExpansion
|Ric_{\alpha\overline{\beta}}-\frac{R}{n}h_{\alpha\overline{\beta}}%
|^{2}\right)  ^{\frac{1}{2}}+\frac{2n}{\lambda_{1}}||\operatorname{Im}%
(A_{\alpha\beta},^{\alpha\beta})||_{L^{2}}%
\end{array}
\label{A1}%
\end{equation}
where $\kappa=\frac{2n\left(  n+1\right)  }{\left(  n-1\right)  \left(
n+2\right)  }(1+\frac{2nK}{\left(  n+1\right)  \lambda_{1}})$ and $\lambda_{1}
$ is the first positive eigenvalue of the sub-Laplacian. Moreover, if the
equality holds then the contact form $\theta$ is pseudo-Einstein,
$\operatorname{Im}(A_{\alpha\beta},^{\alpha\beta})=0,$ and the
pseudo-Hermitian scalar curvature $R=\overline{R}$ is a constant.
\end{theorem}

We observe that when $K=0$ in Theorem \ref{T}, we obtain Theorem \ref{thm}.
Moreover, equality holds in (\ref{1}), we know that the contact form $\theta$
is pseudo-Einstein and $R=\overline{R}$ is a constant. This result is stronger
than we gave in Theorem \ref{thm}.

In Section $3$, we consider a closed pseudo-Hermitian $(2n+1)$-manifold $M$
with zero pseudo-Hermitian torsion and we derive a lower bound estimate for
the first positive eigenvalue $\lambda_{1}$ of the sub-Laplacian $\Delta_{b}$
by using the diameter of $M$ and lower bound of Webster Ricci tensor (see
Proposition \ref{pro}). This estimate is also an independent interesting
result. As a consequence, we have the following Corollary.

\begin{corollary}
\label{cor}Under the same conditions as in the Theorem \ref{T}. We also assume
that $M$ is torsion free, then
\begin{equation}%
\begin{array}
[c]{l}%
\int_{M}(R-\overline{R})^{2}\leq C(Kd^{2})\int_{M}\sum_{\alpha,\beta
}|Ric_{\alpha\overline{\beta}}-\frac{R}{n}h_{\alpha\overline{\beta}}|^{2},
\end{array}
\label{A2}%
\end{equation}
where $d$ is the diameter of $M$ with respect to the Carnot-Carath\'{e}odory
distance and $C(Kd^{2})$ is a constant only depending on $Kd^{2}$. Moreover,
the equality holds if and only if the contact form $\theta$ is pseudo-Einstein.
\end{corollary}

\section{The Proof of Theorem \ref{T}}

In this section, we take the method used in \cite{c1} to prove the inequality
(\ref{A1}). But in the equality case, our proof is more different from \cite{c1}.

We need the following two integral formulas. The first integral is equation
(3.4) in \cite{csw}, for $n\geq2$ and a smooth real-valued function $\varphi
$,
\begin{equation}%
\begin{array}
[c]{ll}
& \frac{n+2}{n-1}\int_{M}\sum_{\alpha,\beta}|\varphi_{\overline{\alpha}\beta
}-\frac{1}{n}\varphi^{\gamma}{}_{\gamma}h_{\overline{\alpha}\beta}|^{2}\\
= & \frac{n+1}{2n}\int_{M}\left(  \Delta_{b}\varphi\right)  ^{2}-\int_{M}%
\sum_{\alpha,\beta}\varphi_{\alpha\beta}\varphi_{\overline{\alpha}%
\overline{\beta}}-\int_{M}(Ric-\frac{n+1}{2}Tor)\left(  \left(  \nabla
_{b}\varphi\right)  _{\mathbb{C}},\left(  \nabla_{b}\varphi\right)
_{\mathbb{C}}\right)  .
\end{array}
\label{5}%
\end{equation}
The second integral comes from Lemma 2.2 in \cite{cc} with its last equation
in P. 268,
\begin{equation}%
\begin{array}
[c]{c}%
\frac{n^{2}}{2}\int_{M}\varphi_{0}^{2}=\int_{M}\sum_{\alpha,\beta}%
(\varphi_{\alpha\beta}\varphi_{\overline{\alpha}\overline{\beta}}%
-\varphi_{\alpha\overline{\beta}}\varphi_{\overline{\alpha}\beta})+\int
_{M}(Ric+\frac{n}{2}Tor)\left(  \left(  \nabla_{b}\varphi\right)
_{\mathbb{C}},\left(  \nabla_{b}\varphi\right)  _{\mathbb{C}}\right)  ,
\end{array}
\label{5a}%
\end{equation}
where $\varphi_{0}=T\varphi$ and $T$ is the characteristic vector field of the
contact form $\theta.$

\emph{Proof of Theorem \ref{T}:}%

%TCIMACRO{\TeXButton{Proof}{\proof}}%
%BeginExpansion
\proof
%EndExpansion
We denote the traceless Webster Ricci tensor by$\ \mathring{R}ic_{\alpha
\overline{\beta}}=Ric_{\alpha\overline{\beta}}-\frac{R}{n}h_{\alpha
\overline{\beta}},$ then the contracted Bianchi identity yields%
\begin{equation}%
\begin{array}
[c]{l}%
\mathring{R}ic_{\alpha\overline{\beta}},^{\overline{\beta}}=\left(
Ric_{\alpha\overline{\beta}}-\frac{R}{n}h_{\alpha\overline{\beta}}\right)
,^{\overline{\beta}}=\frac{n-1}{n}(R_{\alpha}-inA_{\alpha\beta},^{\beta}).
\end{array}
\label{6}%
\end{equation}
Let $f$ be the unique solution of $\Delta_{b}f=R-\overline{R}$ with $\int
_{M}f=0.$ According to (\ref{6}), we compute%

\begin{equation}%
\begin{array}
[c]{ll}
& \int_{M}\left(  R-\overline{R}\right)  ^{2}\\
= & \int_{M}\left(  R-\overline{R}\right)  \Delta_{b}f=-\int_{M}\left\langle
\nabla_{b}R,\nabla_{b}f\right\rangle =-\int_{M}(R_{\alpha}f^{\alpha
}+R_{\overline{\alpha}}f^{\overline{\alpha}})\\
= & \left(  -\frac{n}{n-1}\int_{M}\mathring{R}ic_{\alpha\overline{\beta}%
},^{\overline{\beta}}f^{\alpha}-in\int_{M}A_{\alpha\beta},^{\beta}f^{\alpha
}\right)  +\text{\textrm{complex conjugate}}\\
= & \left(  \frac{n}{n-1}\int_{M}\mathring{R}ic_{\alpha\overline{\beta}%
}f^{\alpha\overline{\beta}}+in\int_{M}A_{\alpha\beta}f^{\alpha\beta}\right)
+\text{\textrm{complex conjugate}}\\
= & \left(  \frac{n}{n-1}\int_{M}\mathring{R}ic_{\alpha\overline{\beta}%
}(f^{\alpha\overline{\beta}}-\frac{1}{n}f^{\gamma}{}_{\gamma}h^{\alpha
\overline{\beta}})+in\int_{M}A_{\alpha\beta}f^{\alpha\beta}\right)
+\text{\textrm{complex conjugate}}\\
= & \frac{2n}{n-1}\int_{M}\mathring{R}ic_{\alpha\overline{\beta}}%
(f^{\alpha\overline{\beta}}-\frac{1}{n}f^{\gamma}{}_{\gamma}h^{\alpha
\overline{\beta}})+in\int_{M}(A_{\alpha\beta}f^{\alpha\beta}-A_{\overline
{\alpha}\overline{\beta}}f^{\overline{\alpha}\overline{\beta}})\\
\leq & \frac{2n}{n-1}||\mathring{R}ic||_{L^{2}}\left(  \int_{M}\sum
_{\alpha,\beta}|f_{\overline{\alpha}\beta}-\frac{1}{n}f^{\gamma}{}_{\gamma
}h_{\overline{\alpha}\beta}|^{2}\right)  ^{\frac{1}{2}}+in\int_{M}%
(A_{\alpha\beta}f^{\alpha\beta}-A_{\overline{\alpha}\overline{\beta}%
}f^{\overline{\alpha}\overline{\beta}}),
\end{array}
\label{20}%
\end{equation}
here we used $\int_{M}\mathring{R}ic_{\alpha\overline{\beta}}f^{\gamma}%
{}_{\gamma}h^{\alpha\overline{\beta}}=0$ and $f^{\alpha\overline{\beta}}%
-\frac{1}{n}f^{\gamma}{}_{\gamma}h^{\alpha\overline{\beta}}=f^{\overline
{\beta}\alpha}-\frac{1}{n}f^{\overline{\gamma}}{}_{\overline{\gamma}%
}h^{\overline{\beta}\alpha}$ is symmetric in $\alpha,\beta$.

Now from (\ref{5}) and the assumption on the curvature condition (\ref{0}), we
obtain
\begin{equation}%
\begin{array}
[c]{l}%
\frac{n+2}{n-1}\int_{M}\sum_{\alpha,\beta}|f_{\overline{\alpha}\beta}-\frac
{1}{n}f^{\gamma}{}_{\gamma}h_{\overline{\alpha}\beta}|^{2}\leq\frac{n+1}%
{2n}\int_{M}\left(  \Delta_{b}f\right)  ^{2}+K\int_{M}\left\vert \nabla
_{b}f\right\vert ^{2}.
\end{array}
\label{21}%
\end{equation}
Besides, by using integration by parts and H\"{o}lder inequality,%
\begin{equation}%
\begin{array}
[c]{c}%
i\int_{M}(A_{\alpha\beta}f^{\alpha\beta}-A_{\overline{\alpha}\overline{\beta}%
}f^{\overline{\alpha}\overline{\beta}})=i\int_{M}(A_{\alpha\beta}%
,^{\alpha\beta}-A_{\overline{\alpha}\overline{\beta}},^{\overline{\alpha
}\overline{\beta}})f\leq2||\operatorname{Im}(A_{\alpha\beta},^{\alpha\beta
})||_{L^{2}}\left\Vert f\right\Vert _{L^{2}}.
\end{array}
\label{21-1}%
\end{equation}
Since the first positive eigenvalue $\lambda_{1}$ of sub-Laplacian on $M$ is
characterized by%
\[%
\begin{array}
[c]{c}%
\lambda_{1}=\inf\left\{  \int_{M}\left\vert \nabla_{b}\varphi\right\vert
^{2}/\int_{M}\varphi^{2}|\text{ }\varphi\text{ \textrm{is nontrivial and}
}\int_{M}\varphi=0\right\}  ,
\end{array}
\]
we have%
\[%
\begin{array}
[c]{lll}%
\int_{M}\left\vert \nabla_{b}f\right\vert ^{2} & = & -\int_{M}f\Delta
_{b}f=-\int_{M}f\left(  R-\overline{R}\right)  \leq\left\Vert f\right\Vert
_{L^{2}}\left\Vert R-\overline{R}\right\Vert _{L^{2}}\\
& \leq & \lambda_{1}^{-1/2}||\nabla_{b}f||_{L^{2}}\left\Vert R-\overline
{R}\right\Vert _{L^{2}}.
\end{array}
\]
Then
\begin{equation}%
\begin{array}
[c]{l}%
\lambda_{1}\int_{M}\left\vert \nabla_{b}f\right\vert ^{2}\leq\left\Vert
R-\overline{R}\right\Vert _{L^{2}}^{2}\text{ \ \textrm{and } }\lambda_{1}%
^{2}\int_{M}\left\vert f\right\vert ^{2}\leq\left\Vert R-\overline
{R}\right\Vert _{L^{2}}^{2}.
\end{array}
\label{22}%
\end{equation}
Due to (\ref{22}), we can rewrite (\ref{21}) and (\ref{21-1}) as
\[%
\begin{array}
[c]{c}%
\frac{n+2}{n-1}\int_{M}\sum_{\alpha,\beta}|f_{\overline{\alpha}\beta}-\frac
{1}{n}f^{\gamma}{}_{\gamma}h_{\overline{\alpha}\beta}|^{2}\leq(\frac{n+1}%
{2n}+\frac{K}{\lambda_{1}})||R-\overline{R}||_{L^{2}}^{2}%
\end{array}
\]
and%
\[%
\begin{array}
[c]{c}%
i\int_{M}(A_{\alpha\beta}f^{\alpha\beta}-A_{\overline{\alpha}\overline{\beta}%
}f^{\overline{\alpha}\overline{\beta}})\leq\frac{2}{\lambda_{1}}%
||\operatorname{Im}(A_{\alpha\beta},^{\alpha\beta})||_{L^{2}}\left\Vert
R-\overline{R}\right\Vert _{L^{2}}%
\end{array}
\]
which combine with (\ref{20}), we then give the equation (\ref{A1}).

Moreover, if the equality of (\ref{A1}) holds, then $f$ will satisfy

\begin{itemize}
\item[(i)] $(Ric-\frac{n+1}{2}Tor+2K)\left(  \left(  \nabla_{b}f\right)
_{\mathbb{C}},\left(  \nabla_{b}f\right)  _{\mathbb{C}}\right)  =0,$

\item[(ii)] $f_{\alpha\beta}=0$ for all $\alpha,\beta,$

\item[(iii)] $R-\overline{R}=c_{1}f$ and $\operatorname{Im}(A_{\alpha\beta
},^{\alpha\beta})=c_{2}f$ for some real constants $c_{1}$ and $c_{2}$,

\item[(iv)] $f_{\alpha\overline{\beta}}-\frac{1}{n}f_{\gamma}$ $^{\gamma
}h_{\alpha\overline{\beta}}=\mu\mathring{R}ic_{\alpha\overline{\beta}}$ for
some constant $\mu$, and

\item[(v)] $\lambda_{1}\int_{M}f^{2}=\int_{M}\left\vert \nabla_{b}f\right\vert
^{2}$ and $\lambda_{1}\int_{M}\left\vert \nabla_{b}f\right\vert ^{2}=\int
_{M}\left(  R-\overline{R}\right)  ^{2}.$
\end{itemize}

Simple computation shows that
\begin{equation}%
\begin{array}
[c]{l}%
\Delta_{b}f=R-\overline{R}=-\lambda_{1}f
\end{array}
\label{30}%
\end{equation}
and
\begin{equation}%
\begin{array}
[c]{c}%
f_{\alpha\overline{\beta}}-\frac{1}{n}f_{\gamma}\text{ }^{\gamma}%
h_{\alpha\overline{\beta}}=\mu\mathring{R}ic_{\alpha\overline{\beta}}%
\end{array}
\label{31}%
\end{equation}
with $\mu=\frac{n+1}{n+2}(1+\frac{2nK}{\left(  n+1\right)  \lambda_{1}}).$

In order to show $\theta$ is pseudo-Einstein, $R=\overline{R}$ is a constant
and $\operatorname{Im}(A_{\alpha\beta},^{\alpha\beta})=0$, it suffices to
claim that $f$ is identically zero. So we need to derive some equations from
(i)$\sim$(v). First, we claim that%

\begin{equation}%
\begin{array}
[c]{c}%
Ric_{\alpha\overline{\beta}}f^{\alpha}-i(n+2)\mu A_{\overline{\alpha}%
\overline{\beta}},^{\overline{\beta}}-i(n+1)A_{\overline{\alpha}%
\overline{\beta}}f^{\overline{\beta}}+2Kf_{\overline{\beta}}=0.
\end{array}
\label{32}%
\end{equation}
We differentiate (ii) and use (\ref{31}), we have%
\[%
\begin{array}
[c]{lll}%
0 & = & f_{\alpha\beta\overline{\gamma}}=f_{\alpha\overline{\gamma}\beta
}+ih_{\beta\overline{\gamma}}f_{\alpha0}+R_{\alpha}{}^{\rho}{}_{\beta
\overline{\gamma}}f_{\rho}\\
& = & \frac{1}{n}f_{\sigma}^{\text{\ \ }\sigma}{}_{\beta}h_{\alpha
\overline{\gamma}}+\mu\mathring{R}ic_{\alpha\overline{\gamma}},_{\beta
}+ih_{\beta\overline{\gamma}}f_{\alpha0}+R_{\alpha}{}^{\rho}{}_{\beta
\overline{\gamma}}f_{\rho}.
\end{array}
\]
Contracting with $h^{\beta\overline{\gamma}},$ we obtain%
\begin{equation}%
\begin{array}
[c]{l}%
0=\frac{1}{n}f_{\sigma}^{\text{ \ }\sigma}{}_{\alpha}+\mu\mathring
{R}ic_{\alpha\overline{\beta}},^{\overline{\beta}}+inf_{\alpha0}%
+Ric_{\alpha\overline{\beta}}f^{\overline{\beta}}.
\end{array}
\label{33}%
\end{equation}
By differentiating the equation (\ref{30}) yields%
\[%
\begin{array}
[c]{lll}%
-\lambda_{1}f_{\alpha} & = & f_{\sigma}^{\text{ \ }\sigma}{}_{\alpha
}+f_{\overline{\sigma}}^{\text{ \ }\overline{\sigma}}{}_{\alpha}=f_{\sigma
}^{\text{ \ }\sigma}{}_{\alpha}+P_{\alpha}f-inA_{\alpha\beta}f^{\beta}\\
& = & f_{\sigma}^{\text{ \ }\sigma}{}_{\alpha}+\frac{n}{n-1}(f_{\alpha
\overline{\beta}}-\frac{1}{n}f_{\gamma}^{\text{ \ }\gamma}h_{\alpha
\overline{\beta}}),^{\overline{\beta}}-inA_{\alpha\beta}f^{\beta}\\
& = & f_{\sigma}^{\text{ \ }\sigma}{}_{\alpha}+\frac{n\mu}{n-1}\mathring
{R}ic_{\alpha\overline{\beta}},^{\overline{\beta}}-inA_{\alpha\beta}f^{\beta},
\end{array}
\]
here the operator $P_{\alpha}f$ is defined by $P_{\alpha}f=f_{\overline
{\sigma}}^{\text{ \ }\overline{\sigma}}{}_{\alpha}+inA_{\alpha\beta}f^{\beta}$
and the second equation follows from equation (3.3) in \cite{gl}. Thus, the
contracted Bianchi identity (\ref{6}) and $R_{\alpha}=-\lambda_{1}f_{\alpha}$
which follows from (\ref{30}) imply
\[%
\begin{array}
[c]{lll}%
f_{\sigma}^{\text{ \ }\sigma}{}_{\alpha} & = & -\lambda_{1}f_{\alpha}%
-\frac{n\mu}{n-1}\mathring{R}ic_{\alpha\overline{\beta}},^{\overline{\beta}%
}+inA_{\alpha\beta}f^{\beta}\\
& = & -\lambda_{1}f_{\alpha}-\mu(R_{\alpha}-inA_{\alpha\beta},^{\beta
})+inA_{\alpha\beta}f^{\beta}\\
& = & (\mu-1)\lambda_{1}f_{\alpha}+in\mu A_{\alpha\beta},^{\beta}%
+inA_{\alpha\beta}f^{\beta}.
\end{array}
\]
Also, by the commutation relations (\cite[Lemma 2.3]{le}), we have%
\[%
\begin{array}
[c]{lll}%
inf_{\alpha0} & = & inf_{0\alpha}-inA_{\alpha\beta}f^{\beta}=(f_{\sigma
}^{\text{ \ }\sigma}-f_{\overline{\sigma}}^{\text{ \ }\overline{\sigma}%
}),_{\alpha}-inA_{\alpha\beta}f^{\beta}\\
& = & f_{\sigma}^{\text{ \ }\sigma}{}_{\alpha}-P_{\alpha}f=f_{\sigma}^{\text{
\ }\sigma}{}_{\alpha}-\frac{n\mu}{n-1}\mathring{R}ic_{\alpha\overline{\beta}%
},^{\overline{\beta}}\\
& = & (2\mu-1)\lambda_{1}f_{\alpha}+2in\mu A_{\alpha\beta},^{\beta
}+inA_{\alpha\beta}f^{\beta}.
\end{array}
\]
Substituting these into (\ref{33}) and using the fact $\mu=\frac{n+1}%
{n+2}(1+\frac{2nK}{\left(  n+1\right)  \lambda_{1}})$, we final get%
\[%
\begin{array}
[c]{lll}%
0 & = & \frac{1}{n}f_{\sigma}^{\text{ \ }\sigma}{}_{\alpha}+\mu\mathring
{R}ic_{\alpha\overline{\beta}},^{\overline{\beta}}+inf_{\alpha0}%
+Ric_{\alpha\overline{\beta}}f^{\overline{\beta}}\\
& = & 2Kf_{\alpha}+i(n+2)\mu A_{\alpha\beta},^{\beta}+i(n+1)A_{\alpha\beta
}f^{\beta}+Ric_{\alpha\overline{\beta}}f^{\overline{\beta}},
\end{array}
\]
which is (\ref{32}) as claimed.

Next, we want to show
\begin{equation}%
\begin{array}
[c]{c}%
\int_{M}A_{\alpha\beta}f^{\alpha}f^{\beta}=\int_{M}A_{\overline{\alpha
}\overline{\beta}}f^{\overline{\alpha}}f^{\overline{\beta}}=0\text{
\ }\mathrm{and}\text{ \ }\int_{M}Ric_{\alpha\overline{\beta}}f^{\alpha
}f^{\overline{\beta}}+K\int_{M}\left\vert \nabla_{b}f\right\vert ^{2}=0.
\end{array}
\label{41}%
\end{equation}
From (\ref{32}) we know that%
\[%
\begin{array}
[c]{l}%
Ric_{\alpha\overline{\beta}}f^{\alpha}f^{\overline{\beta}}+i\left(
n+2\right)  \mu A_{\alpha\beta},^{\beta}f^{\alpha}+i(n+1)A_{\alpha\beta
}f^{\alpha}f^{\beta}+2Kf_{\alpha}f^{\alpha}=0.
\end{array}
\]
But compare this with (i)
\begin{equation}%
\begin{array}
[c]{c}%
Ric_{\alpha\overline{\beta}}f^{\alpha}f^{\overline{\beta}}+\frac{n+1}%
{2}i(A_{\alpha\beta}f^{\alpha}f^{\beta}-A_{\overline{\alpha}\overline{\beta}%
}f^{\overline{\alpha}}f^{\overline{\beta}})+2Kf_{\alpha}f^{\alpha}=0,
\end{array}
\label{42}%
\end{equation}
one gets
\[%
\begin{array}
[c]{c}%
\frac{2(n+2)}{n+1}\mu A_{\alpha\beta},^{\beta}f^{\alpha}=-(A_{\alpha\beta
}f^{\alpha}f^{\beta}+A_{\overline{\alpha}\overline{\beta}}f^{\overline{\alpha
}}f^{\overline{\beta}}).
\end{array}
\]
Then integral it yields
\[%
\begin{array}
[c]{c}%
\int_{M}A_{\alpha\beta}f^{\alpha}f^{\beta}+\int_{M}A_{\overline{\alpha
}\overline{\beta}}f^{\overline{\alpha}}f^{\overline{\beta}}=0
\end{array}
\]
due to $\int_{M}A_{\alpha\beta},^{\beta}f^{\alpha}=\int_{M}(A_{\alpha\beta
}f^{\alpha}),^{\beta}=0,$ by (ii). Also by the reality of $A_{\alpha\beta
},^{\beta}f^{\alpha}$, we know
\[%
\begin{array}
[c]{c}%
\int_{M}A_{\alpha\beta}f^{\alpha}f^{\beta}=-\int_{M}A_{\overline{\alpha
}\overline{\beta}}f^{\overline{\alpha}}f^{\overline{\beta}}=\int
_{M}A_{\overline{\alpha}\overline{\beta}},^{\overline{\beta}}f^{\overline
{\alpha}}f
\end{array}
\]
is real. Hence, the integral of (\ref{42}),
\[%
\begin{array}
[c]{c}%
\int_{M}Ric_{\alpha\overline{\beta}}f^{\alpha}f^{\overline{\beta}}+\left(
n+1\right)  i\int_{M}A_{\alpha\beta}f^{\alpha}f^{\beta}+K\int_{M}\left\vert
\nabla_{b}f\right\vert ^{2}=0
\end{array}
\]
will imply (\ref{41}) as we wanted.

Now, by applying (ii) and (\ref{41}) to the equation (\ref{5a}), we final
obtain
\[%
\begin{array}
[c]{c}%
\frac{n^{2}}{2}\int_{M}f_{0}^{2}+\int_{M}\sum_{\alpha,\beta}f_{\alpha
\overline{\beta}}f_{\overline{\alpha}\beta}+K\int_{M}\left\vert \nabla
_{b}f\right\vert ^{2}=0.
\end{array}
\]
It implies that $f=0$ as desired. This completes the proof of Theorem
\ref{T}.
%TCIMACRO{\TeXButton{End Proof}{\endproof}}%
%BeginExpansion
\endproof
%EndExpansion

\section{First eigenvalue estimate of the sub-Laplacian}

Let $\left(  M,J,\theta\right)  $ be a closed pseudo-Hermitian $(2n+1)$%
-manifold with vanishing pseudo-Hermitian torsion. In this section, by
applying the argument of the CR analogous Li-Yau's gradient estimate in
\cite{ckl}, we derive a lower bound estimate for the first positive eigenvalue
$\lambda_{1}$ of the sub-Laplacian $\Delta_{b}$ using the diameter of $M$ and
lower bound of Webster Ricci tensor. In the last, we prove Corollary \ref{cor}.

\begin{proposition}
\label{pro} Let $\left(  M,J,\theta\right)  $ be a closed pseudo-Hermitian
$(2n+1)$-manifold with vanishing pseudo-Hermitian torsion and the Webster
Ricci tensor is bounded from below by a nonpositive constant $-K.$ Let $f$ be
an eigenfunction of $\Delta_{b}$ with respect to the first eigenvalue
$\lambda_{1}.$ Then there exists constants $C_{1}(n),$ $C_{2}(n)>0$ depending
on $n$ alone, such that%
\begin{equation}%
\begin{array}
[c]{c}%
\lambda_{1}\geq\frac{C_{1}}{d^{2}}\exp(-C_{2}d\sqrt{K}).
\end{array}
\label{B}%
\end{equation}
Here $d$ is the diameter of $M$ with respect to the Carnot-Carath\'{e}odory
distance (see definition 2.3 in \cite{ckl}).
\end{proposition}

We recall the following CR version of Bochner formula in a pseudo-Hermitian
$(2n+1)$-manifold (\cite{gr}). For a smooth real-valued function $\varphi,$
\begin{equation}%
\begin{array}
[c]{lll}%
\frac{1}{2}\Delta_{b}\left\vert \nabla_{b}\varphi\right\vert ^{2} & = &
|(\nabla^{H})^{2}\varphi|^{2}+\left\langle \nabla_{b}\varphi,\nabla_{b}%
\Delta_{b}\varphi\right\rangle +2\left\langle J\nabla_{b}\varphi,\nabla
_{b}\varphi_{0}\right\rangle \\
&  & +[2Ric-(n-2)Tor]\left(  \left(  \nabla_{b}\varphi\right)  _{\mathbb{C}%
},\left(  \nabla_{b}\varphi\right)  _{\mathbb{C}}\right)  .
\end{array}
\label{00}%
\end{equation}
Since
\[%
\begin{array}
[c]{c}%
|(\nabla^{H})^{2}\varphi|^{2}=2\sum_{\alpha,\beta}(|\varphi_{\alpha\beta}%
|^{2}+|\varphi_{\alpha\overline{\beta}}|^{2})\geq2\sum_{\alpha}|\varphi
_{\alpha\overline{\alpha}}|^{2}\geq\frac{1}{2n}\left(  \Delta_{b}%
\varphi\right)  ^{2}+\frac{n}{2}\varphi_{0}^{2}%
\end{array}
\]
and for any constant $v>0$,%
\[
2\left\langle J\nabla_{b}\varphi,\nabla_{b}\varphi_{0}\right\rangle
\leq2\left\vert \nabla_{b}\varphi\right\vert \left\vert \nabla_{b}\varphi
_{0}\right\vert \leq v^{-1}\left\vert \nabla_{b}\varphi\right\vert
^{2}+v\left\vert \nabla_{b}\varphi_{0}\right\vert ^{2}.
\]
Therefore, for a real function $\varphi$ and any $v>0$, the Bochner formula
(\ref{00}) becomes
\begin{equation}%
\begin{array}
[c]{lll}%
\Delta_{b}\left\vert \nabla_{b}\varphi\right\vert ^{2} & \geq & \frac{1}%
{n}\left(  \Delta_{b}\varphi\right)  ^{2}+n\varphi_{0}^{2}+2\left\langle
\nabla_{b}\varphi,\nabla_{b}\Delta_{b}\varphi\right\rangle -2v\left\vert
\nabla_{b}\varphi_{0}\right\vert ^{2}\\
&  & +2[2Ric-(n-2)Tor-2v^{-1}]\left(  \left(  \nabla_{b}\varphi\right)
_{\mathbb{C}},\left(  \nabla_{b}\varphi\right)  _{\mathbb{C}}\right)  .
\end{array}
\label{2}%
\end{equation}

\emph{Proof of Proposition} \ref{pro}:%

%TCIMACRO{\TeXButton{Proof}{\proof} }%
%BeginExpansion
\proof
%EndExpansion
Let $f$ be an eigenfunction of $\Delta_{b}$ with respect to the eigenvalue
$\lambda_{1}$. Since
\[%
\begin{array}
[c]{c}%
\lambda_{1}\int_{M}f=-\int_{M}\Delta_{b}f=0,
\end{array}
\]
$f$ must change sign. We may normalize $f$ to satisfy $\min f=-1$ and $\max
f\leq1.$ Let us consider the function $\varphi=\ln(f+a),$ for some constant
$a>1.$ Then the function $\varphi$ satisfies%
\[%
\begin{array}
[c]{c}%
\Delta_{b}\varphi=-\left\vert \nabla_{b}\varphi\right\vert ^{2}-\frac
{\lambda_{1}f}{f+a}%
\end{array}
\]
and thus
\begin{equation}%
\begin{array}
[c]{l}%
\left\langle \nabla_{b}\varphi,\nabla_{b}\Delta_{b}\varphi\right\rangle
=-\left\langle \nabla_{b}\varphi,\nabla_{b}\left\vert \nabla_{b}%
\varphi\right\vert ^{2}\right\rangle -\frac{a\lambda_{1}}{f+a}\left\vert
\nabla_{b}\varphi\right\vert ^{2}.
\end{array}
\label{11}%
\end{equation}
Since
\[%
\begin{array}
[c]{lll}%
\Delta_{b}\varphi_{0} & = & \left(  \Delta_{b}\varphi\right)  _{0}%
+2[(A_{\alpha\beta}\varphi^{\beta}),^{\alpha}+(A_{\overline{\alpha}%
\overline{\beta}}\varphi^{\overline{\beta}}),^{\overline{\alpha}}]\\
& = & \left(  -\left\vert \nabla_{b}\varphi\right\vert ^{2}-\frac{\lambda
_{1}f}{f+a}\right)  _{0}=-2\left\langle \nabla_{b}\varphi,\nabla_{b}%
\varphi_{0}\right\rangle -\frac{a\lambda_{1}}{f+a}\varphi_{0}.
\end{array}
\]
Therefore, we have
\begin{equation}%
\begin{array}
[c]{c}%
\frac{1}{2}\Delta_{b}\varphi_{0}^{2}=\left\vert \nabla_{b}\varphi
_{0}\right\vert ^{2}+\varphi_{0}\Delta_{b}\varphi_{0}=\left\vert \nabla
_{b}\varphi_{0}\right\vert ^{2}-\left\langle \nabla_{b}\varphi,\nabla
_{b}\varphi_{0}^{2}\right\rangle -\frac{a\lambda_{1}}{f+a}\varphi_{0}^{2}.
\end{array}
\label{12}%
\end{equation}
And
\begin{equation}%
\begin{array}
[c]{c}%
\Delta_{b}\frac{f}{f+a}=\frac{a}{\left(  f+a\right)  ^{2}}\Delta_{b}%
f-\frac{2a}{\left(  f+a\right)  ^{3}}\left\vert \nabla_{b}f\right\vert
^{2}=-2\left\langle \nabla_{b}\varphi,\nabla_{b}\frac{f}{f+a}\right\rangle
-\frac{a\lambda_{1}}{f+a}\frac{f}{f+a}.
\end{array}
\label{10}%
\end{equation}

Now, we define $F:M\times\lbrack0,1]\rightarrow\mathbb{R}$ by
\[%
\begin{array}
[c]{lll}%
F\left(  x,t\right)  & = & t(\left\vert \nabla_{b}\varphi\right\vert
^{2}-\alpha\frac{\lambda_{1}f}{f+a}+\gamma t\varphi_{0}^{2})\\
& = & t[(\alpha+1)\left\vert \nabla_{b}\varphi\right\vert ^{2}+\alpha
\Delta_{b}\varphi+\gamma t\varphi_{0}^{2}],
\end{array}
\]
where $\alpha$ be a nonzero constant and $\gamma$\ be a positive constant
which will be chosen later. By applying the Bochner inequality (\ref{2}) with
$v=\gamma t$, and using (\ref{11}), (\ref{12}) and (\ref{10}), one can derive%
\begin{equation}%
\begin{array}
[c]{lll}%
\Delta_{b}F+2\left\langle \nabla_{b}\varphi,\nabla_{b}F\right\rangle  & \geq &
\frac{t}{n}\left(  \Delta_{b}\varphi\right)  ^{2}+nt\varphi_{0}^{2}-2\left(
Kt+\gamma^{-1}\right)  \left\vert \nabla_{b}\varphi\right\vert ^{2}\\
&  & -\frac{a\lambda_{1}}{f+a}t[2\left\vert \nabla_{b}\varphi\right\vert
^{2}-\alpha\frac{\lambda_{1}f}{f+a}+2\gamma t\varphi_{0}^{2}]\\
& = & \frac{t}{n}\left(  \Delta_{b}\varphi\right)  ^{2}+nt\varphi_{0}%
^{2}-2\left(  Kt+\gamma^{-1}\right)  \left\vert \nabla_{b}\varphi\right\vert
^{2}\\
&  & -\frac{a\lambda_{1}}{f+a}[F+t\left\vert \nabla_{b}\varphi\right\vert
^{2}+\gamma t^{2}\varphi_{0}^{2}].
\end{array}
\label{13}%
\end{equation}

On the other hand, from the definition of $F(x,t)$, we have
\[%
\begin{array}
[c]{c}%
\Delta_{b}\varphi=\alpha^{-1}[t^{-1}F-(\alpha+1)\left\vert \nabla_{b}%
\varphi\right\vert ^{2}-\gamma t\varphi_{0}^{2}],
\end{array}
\]
thus%
\[%
\begin{array}
[c]{l}%
\left(  \Delta_{b}\varphi\right)  ^{2}\geq(\alpha t)^{-2}F^{2}-2\alpha
^{-2}t^{-1}F[(\alpha+1)\left\vert \nabla_{b}\varphi\right\vert ^{2}+\gamma
t\varphi_{0}^{2}].
\end{array}
\]
Substituting this into (\ref{13}), we obtain%
\[%
\begin{array}
[c]{ll}
& \Delta_{b}F+2\left\langle \nabla_{b}\varphi,\nabla_{b}F\right\rangle \\
\geq & \frac{1}{n\alpha^{2}t}F^{2}-\frac{a\lambda_{1}}{f+a}F+[n-(\frac
{a\lambda_{1}}{f+a}t+\frac{2}{n\alpha^{2}}F)\gamma]t\varphi_{0}^{2}\\
& -2[Kt+\gamma^{-1}+\frac{t}{2}\frac{a\lambda_{1}}{f+a}+\frac{\alpha
+1}{n\alpha^{2}}F]\left\vert \nabla_{b}\varphi\right\vert ^{2}.
\end{array}
\]
Thus, at a maximum point $p_{t}$ of $F(\cdot,t)$, we have
\[
0\geq\Delta_{b}F(p_{t},t)+2\left\langle \nabla_{b}\varphi,\nabla
_{b}F\right\rangle (p_{t},t).
\]
Hence, at $(p_{t},t),$
\begin{equation}%
\begin{array}
[c]{lll}%
0 & \geq & \frac{1}{n\alpha^{2}t}F^{2}-\frac{a\lambda_{1}}{f+a}F+[n-(\frac
{a\lambda_{1}}{f+a}t+\frac{2}{n\alpha^{2}}F)\gamma]t\varphi_{0}^{2}\\
&  & -2[Kt+\gamma^{-1}+\frac{t}{2}\frac{a\lambda_{1}}{f+a}+\frac{\alpha
+1}{n\alpha^{2}}F]\left\vert \nabla_{b}\varphi\right\vert ^{2}.
\end{array}
\label{14}%
\end{equation}

We claim that there exist constants $\alpha$ depending only on $n$ with
$(\alpha+1)<0$ and $\gamma$ depending on $\lambda_{1}$, $a$ and $K$ such that
\[%
\begin{array}
[c]{c}%
F\left(  x,t\right)  <-\frac{n\alpha^{2}}{\alpha+1}[K+\gamma^{-1}+\frac{1}%
{2}\frac{a\lambda_{1}}{a-1}]
\end{array}
\]
on $M\times\lbrack0,1].$ We prove it by contradiction. Suppose not, then
\[%
\begin{array}
[c]{c}%
\max_{M\times\lbrack0,1]}F(x,t)\geq-\frac{n\alpha^{2}}{\alpha+1}[K+\gamma
^{-1}+\frac{1}{2}\frac{a\lambda_{1}}{a-1}].
\end{array}
\]
Since $F$ is continuous in the variable $t$ and $F\left(  x,0\right)  =0,$
thus there exists a $t_{0}\in(0,1]$ such that
\[%
\begin{array}
[c]{c}%
\max_{M\times\lbrack0,t_{0}]}F(x,t)=-\frac{n\alpha^{2}}{\alpha+1}%
[K+\gamma^{-1}+\frac{1}{2}\frac{a\lambda_{1}}{a-1}].
\end{array}
\]
Assume $F$ achieves its maximum at the point $(p_{t_{0}},s_{t_{0}})$ on
$M\times\lbrack0,t_{0}].$ Then%
\begin{equation}%
\begin{array}
[c]{c}%
F(p_{t_{0}},s_{t_{0}})=-\frac{n\alpha^{2}}{\alpha+1}[K+\gamma^{-1}+\frac{1}%
{2}\frac{a\lambda_{1}}{a-1}]>0.
\end{array}
\label{14'}%
\end{equation}
By applying (\ref{14}) at a maximum point $p_{t_{0}}$ of $F(\cdot,s_{t_{0}}) $
and using (\ref{14'}), one obtain%
\begin{equation}%
\begin{array}
[c]{lll}%
0 & \geq & \frac{1}{n\alpha^{2}s_{t_{0}}}F(p_{t_{0}},s_{t_{0}})^{2}%
-\frac{a\lambda_{1}}{f+a}F(p_{t_{0}},s_{t_{0}})\\
&  & +[n-(\frac{a\lambda_{1}}{f+a}s_{t_{0}}+\frac{2}{n\alpha^{2}}F(p_{t_{0}%
},s_{t_{0}}))\gamma]s_{t_{0}}\varphi_{0}^{2}.
\end{array}
\label{15}%
\end{equation}

Now we choose
\[%
\begin{array}
[c]{c}%
\alpha+1=-\frac{3}{n}\text{ \ }\mathrm{and}\text{ \ }\gamma^{-1}=\frac{n+3}%
{n}\frac{a\lambda_{1}}{a-1}+2K,
\end{array}
\]
then%
\[%
\begin{array}
[c]{c}%
\frac{1}{n\alpha^{2}s_{t_{0}}}F(p_{t_{0}},s_{t_{0}})^{2}-\frac{a\lambda_{1}%
}{f+a}F(p_{t_{0}},s_{t_{0}})\geq ns_{t_{0}}^{-1}(K+\frac{1}{2}\frac
{a\lambda_{1}}{a-1})F(p_{t_{0}},s_{t_{0}})>0
\end{array}
\]
and
\[%
\begin{array}
[c]{c}%
n-(\frac{a\lambda_{1}}{f+a}s_{t_{0}}+\frac{2}{n\alpha^{2}}F(p_{t_{0}}%
,s_{t_{0}}))\gamma\geq\frac{n}{3}[1-(\frac{n+3}{n}\frac{a\lambda_{1}}%
{a-1}+2K)\gamma]=0.
\end{array}
\]
This leads to a contradiction with (\ref{15}). Therefore, we obtain that%
\[%
\begin{array}
[c]{c}%
F\left(  x,t\right)  <(n+3)^{2}(\frac{n+2}{2n}\frac{a\lambda_{1}}{a-1}+K)
\end{array}
\]
on $M\times\lbrack0,1]$. In particular, at $t=1,$ we have
\[%
\begin{array}
[c]{c}%
\left\vert \nabla_{b}\varphi\right\vert ^{2}+\frac{n+3}{n}\frac{\lambda_{1}%
f}{f+a}+(\frac{n+3}{n}\frac{a\lambda_{1}}{a-1}+2K)^{-1}\varphi_{0}^{2}%
\leq(n+3)^{2}(\frac{n+2}{2n}\frac{a\lambda_{1}}{a-1}+K).
\end{array}
\]
Thus, we obtain the subgradient estimate
\[%
\begin{array}
[c]{lll}%
\left\vert \nabla_{b}\varphi\right\vert ^{2}+(\frac{n+3}{n}\frac{a\lambda_{1}%
}{a-1}+2K)^{-1}\varphi_{0}^{2} & \leq & (n+3)^{2}(\frac{n+2}{2n}\frac
{a\lambda_{1}}{a-1}+K)-\frac{n+3}{n}\frac{\lambda_{1}f}{f+a}\\
& \leq & (n+3)^{2}(\frac{n+3}{2n}\frac{a\lambda_{1}}{a-1}+K).
\end{array}
\]
Therefore,
\[%
\begin{array}
[c]{c}%
\left\vert \nabla_{b}\varphi\right\vert ^{2}\leq(n+3)^{2}(\frac{n+3}{2n}%
\frac{a\lambda_{1}}{a-1}+K).
\end{array}
\]
By integrating $\left\vert \nabla_{b}\varphi\right\vert =\left\vert \nabla
_{b}\ln(f+a)\right\vert $ along a minimal horizontal geodesic $\varsigma$
joining the points at which $f=-1$ and $f=\max f,$ it follows that%
\[%
\begin{array}
[c]{lll}%
\ln\frac{a}{a-1} & \leq & \ln\left(  \frac{a+\max f}{a-1}\right)  =\ln(a+\max
f)-\ln(a-1)\\
& \leq & \int_{\varsigma}\left\vert \nabla_{b}\ln(f+a)\right\vert
\leq(n+3)d\sqrt{\frac{n+3}{2n}\frac{a\lambda_{1}}{a-1}+K},
\end{array}
\]
for all $a>1.$ Setting $s=(a-1)/a,$ we obtain%
\[%
\begin{array}
[c]{c}%
\frac{(n+3)^{3}}{2n}\lambda_{1}\geq\left(  d^{-2}(\ln s^{-1})^{2}%
-(n+3)^{2}K\right)  s
\end{array}
\]
for all $0<s<1.$ Maximizing the right hand side as a function of $s$ by
setting $s=\exp(-1-\sqrt{1+(n+3)^{2}Kd^{2}}),$ we get the estimate
\[%
\begin{array}
[c]{l}%
\lambda_{1}\geq\frac{4n}{(n+3)^{3}d^{2}}\exp(-1-\sqrt{1+(n+3)^{2}Kd^{2}})
\end{array}
\]
as claimed. This completes the proof of Proposition \ref{pro}.
%TCIMACRO{\TeXButton{End Proof}{\endproof}}%
%BeginExpansion
\endproof
%EndExpansion

Now we prove Corollary \ref{cor}.

\emph{Proof} \emph{of Corollary} \ref{cor}: From the Proposition \ref{pro}, we
have%
\[%
\begin{array}
[c]{c}%
\frac{K}{\lambda_{1}}\leq C_{1}Kd^{2}\exp(C_{2}d\sqrt{K}).
\end{array}
\]
So we obtain the inequality in Corollary \ref{cor} with the constant%
\[%
\begin{array}
[c]{c}%
C(Kd^{2})=\frac{4n^{2}}{\left(  n-1\right)  \left(  n+2\right)  }\left(
\frac{n+1}{2n}+C_{1}Kd^{2}\exp(C_{2}d\sqrt{K})\right)  .
\end{array}
\]%
%TCIMACRO{\TeXButton{End Proof}{\endproof}}%
%BeginExpansion
\endproof
%EndExpansion

\subsection*{Acknowledgment}

The second author would like to express his thanks to Prof. Shu-Cheng Chang
for his encouragement and teaching on CR geometry. He is also grateful to
Prof. Xu Cheng for her explanation on the paper \cite{c2}.

\end{document}